\documentclass[11pt]{article}
\usepackage{graphicx}
\usepackage{amssymb}
\usepackage{hyperref}
\usepackage{epstopdf}
\newtheorem{theorem}{Theorem}

\newtheorem{openproblem}{Open Problem}

\def\semi{\hbox{ $\times $ \kern-.972em \raise.12719em\hbox{ $_{^|}$}  }}

\def\bi{\begin{itemize}}
\def\ei{\end{itemize}}
\date{September 9, 2023}
\title{Algebra, Topology and the discoveries of Vaughan Jones  \footnote  {primary classification 57K14 }}
\author{Joan S. Birman}
\begin{document}
\maketitle
\centerline{Dedicated to the memory of Vaughan Jones}

\begin{abstract}  
\noindent In this article  the discovery of the Jones Polynomial will be discussed,  emphasizing the way in which it illustrated the remarkable unity between distinct parts of Mathematics,  each with its own language, but initially without a dictionary.   
\end{abstract}

\section{Introduction}\label{s:introduction}
 In this brief tribute to Vaughan Jones, I will emphasize the way in which the discovery of the Jones polynomial and its generalizations  illustrated the remarkable unity between distinct parts of mathematics,  each with its own language, but initially without a dictionary.  The first language of the dictionary comes from algebra, or more explicitly the subarea of von Neumann algebras  that includes Type II$_1$ factors, while the second part of the dictionary comes from topology, or more expicitly  the subarea that includes knots and links in 3-space.   The dictionary was expanded, through the remarkably rich consequences of Jones' work,  moving outward within both algebra and topology through multiple new discoveries in both.   

My references will include copies of ten  handwritten letters from Jones to me, and one from Jones to Louis Kauffman, that are posted on the web so that interested readers can look at them and witness a piece of the creative process itself.   The letters have not been discussed elsewhere, and   since this article is a memorial tribute to Jones, it seemed like an appropriate place to fill that gap.   

{\bf Acknowledgements:} I thank Louis Kauffman for inviting me to write this paper;  Sofia Lambropoulou for help with the manuscript; Robion Kirby and Sheila Newberry for posting the letters from Jones to me on the {\it Celebratio Mathematica} website; Martha Jones for pointing out a place where my memory regarding the letters was clearly incorrect; and Kenneth Birman, for reading an early version of this manuscript and helping me to improve it.

\section{A knot or link $\bf K $ in $\mathbb R^3$ or $\mathbb S^3$} \label{s:knot polynomials}
My initial focus will be on the material in the papers \cite{Jones-Bull,Jones-Annals} that relate to the one variable polynomial  invariant $V_{\bf K }(t)$ of the type of a classical knot or link $\bf K $ in 3-space, as Vaughan Jones discovered it and discussed it with me.  My account in this section will be a very personal account. The discovery that I will describe in detail ultimately lead to an explosion in knot and link theory, as the part of topology that had centered on Alexander invariants grew into an area that is now known as {\it Quantum Topology}, reaching into both algebra and topology, and resulting in the awarding of a Fields Medal to Jones in 1990.    

My starting point is Jones' foundational paper \cite{Jones-Index}, about towers of type II$_1$ algebras, ordered by inclusion.   See \cite{Bisch et al} for an account of the particular way in which representations of  Artin's braid groups $B_n, n=1,2,3,\dots$  appeared in the work of Jones.  His starting point was a family of irreducible matrices,  representations of the braid groups, and his discovery of a class invariant which he called a `trace'.  It was a weighted sum of traces on the defining matrices.   He had read Artin's paper \cite{Artin}, with its pictures of braids, and was greatly puzzled by the appearance of braids in his work, because the work that he had done on type II$_1$ factors seemed to be so very far from topology.  

Jones  persisted in efforts to understand why they appeared.  In 1983 he gave a talk at a conference in Kyoto, Japan, 
where he described the structure of his representations of Artin's braid groups $B_n, n\geq 1$, The report \cite {Jones-Pitman} that he wrote on his talk contains almost all of the essential building blocks that were needed for the development of the Jones polynomial of a knot or link, in $\mathbb R^3$ or $\mathbb S^3$, however I don't see any indication in it that he realized that his invariant might be related to  the knots and links that one obtains after changing a braid to a `closed braid'.  He persisted in asking mathematicians he met, in several different places, about the possible meaning of representations of $B_n$  in his work on type II$_1$ factors, but (as he told us) met up with dead ends.  

In this regard we note that Geometric Group Theory as we know it today did not exist in 1983, when Vaughan Jones was asking all those questions about his Hecke algebra representations of $B_n$.  This was the case even though  (i) Alexander had proved in 1923  that every knot and link can be represented as a closed braid \cite{Alex 1923}; also (ii) the Alexander polynomial and related structure were at the center of knot theory in 1984, so that Alexander's name and the paper \cite{Alex 1928} were very well known in 1983; also 
(iii) Markov's  1936 Theorem relating any two closed braid representatives of the same link \cite{Markov 1936}, with its light sketch of a proof,  had finally, in 1974,  been backed up fully with a detailed proof in \cite{Birman-book} that was subsequently backed up by other proofs  in (for example) \cite{Morton1986, Lam-Rourke 1997, Traczyk1998}.    

In the early spring of 1984, Jones (on the faculty of U. Penn in Philadelphia) attended a lecture at the nearby Institute for Advanced Study, in Princeton, New Jersey.  Chatting with Caroline Series, then a member, at tea, Jones told Series about his unexpected discovery of representations of braid groups, and she said something like `then you should consult Joan Birman', and that is how we arranged to meet at Columbia University on May 14, 1984.  The date is precise, because it was recorded in a daily calendar that I kept, where I noted two meetings with Jones, the second on May 22, 1984.  Vaughan Jones' discovery of the one-variable Jones polynomial can be pinpointed precisely as having occurred between those two meetings.    

 When he arrived for the May 14 meeting, he told me about his matrix representations of Artin's braid group $B_n$ and his trace. I had written a book about braids \cite{Birman-book}, and Jones had even referenced it in  \cite{Jones-Pitman}.  One of my goals, when I wrote the book, had been to study knots and links by thinking of them as closed braids, and my book had been based on lecture notes from a graduate course that I had given, where I  laid out the foundations for such study.  At our meeting  I pointed out to him the significance of Markov's theorem for that project.  I noted that Markov had stated his theorem, and sketched a proof,  but never published a detailed proof,  however I had filled that gap with a proof in \cite{Birman-book}. 
 
 After our discussions that day he realized that his trace might be related to the representation of classical knots and links in 3-space {\it as closed braids}, rather than to group-theoretic invariants of the braids themselves.  
During the interval between our two meetings Jones  had the key idea that by a simple rescaling of his algebra, he could arrange that his representations satisfied the rules of Markov's Theorem, so that his rescaled trace determined  invariants of the knot or link type $\bf K $.   When he arrived at the May 22 meeting, he had in hand a polynomial invariant of the link type $\bf K $, but he did not know whether it was new.  At the time there was a very well-known and well-studied polynomial invariant of link type, the Alexander polynomial, and in many ways his invariant resembled the Alexander polynomial.     Fortuitously, I had in my office filing cabinet braid representatives that I had computed (out of curiosity) for an infinite  family of braids whose closures had trivial Alexander polynomials.  Jones and I established, together, that Jones' polynomial could not be the Alexander polynomial, because it detected links that had trivial Alexander polynomial.   The Jones polynomial was born.   We both understood that it was bound to be a major discovery in knot theory.  

The discovery of the Jones polynomial was based entirely on algebra.   We pose the following interesting problem: 
\begin{openproblem} Identify topological aspects of knots and/or  links that are captured by the Jones polynomial. 
\end{openproblem}
For example, it is known that the variable $t$ in the Jones polynomial changes to $t^{-1}$ when a knot is replaced by its mirror image, so in this sense {\it amphichirality} of a knot is captured by the Jones polynomial.

\section{The letters}\label{s:letters}

I had plans to leave for a month in England, a few days after the discovery of the one-variable Jones polynomial, later dubbed $V_{\bf K}(t)$ , and  knew that when I returned Jones would be on his way to California for a sabbatical year that he would spend at MSRI, in Berkeley.  The story continued in letters and telephone conversations.  In particular, the reader is invited to go to the website

$  \href{https://celebratio.org/Birman_JS/article/639/}{https://celebratio.org/Birman\_JS/article/639/}  $

where there are handwritten copies of letters to us from Jones. The first letter is dated May 31,1984, a week after our second meeting, and the second letter continued it, and while it is undated, it can easily be seen to have been written the following day. It reviewed (for my benefit) the derivation of what is now known to be the one-variable or $V$-polynomial $V_{\bf K}(t)$ of a knot or link $K\subset \mathbb R^3$ in 3-space. He worked out some of its properties.   
In the fall of 1984 I  resumed my normal academic responsibilities, and was thoroughly involved in and committed to other research projects.   The letters that I found in my office are all one-way, from Vaughan to me.  I recall talking with him  often on the telephone, and also attending a week-long conference at MSRI October 11-16,1984.  (Again, notes in my daily calendar pinpointed the dates).  

Letters like the ones that are posted on the Celebratio Mathematica website offer a look at the creative process as it is occurring.  The first two letters speak of a formula for knots that are closed 3-braids, and  his speculation whether the Jones polynomial might be enough to distinguish them?  That possibility must have  interested me, because during August of 1984 I investigated the matter, and wrote a paper \cite{Birman-3-braids} about it, and learned that in fact, while  $V_{\bf K}(t)$ did a better job than the Alexander polynomial $A_{\bf K}(t)$, it wasn't qualitatively different from $A_{\bf K}(t)$ in that regard.   To this day it is not known whether the Jones polynomial is a complete invariant for a much simpler object than the set of all links of braid index 3, namely an unknotted circle.    

Several of the letters mention  joint work, in passing,  but we never wrote a joint paper.  I was very pleased when I was asked to give an introductory talk on the work of Vaughan Jones  at ICM 1990 in Kyoto, and knew then that Vaughan had acknowledged the role that I had played in his discoveries with appropriate generosity.  

The final letter that's posted  on the website was dated 12 June, 1990.  In the spirit of changing times, it was an e-mail letter rather than a handwritten one.  It addresses a question I had asked my friend Vaughan, regarding how he found the patience to do so many very long calculations without making mistakes?  The reader is invited to look at his 12 June, 1990 letter, to learn his answer about one such calculation.

\section {A space of all knots in 3-space} \label{s:Vassiliev invariants}  We change the focus, and move to a new application of the Jones polynomial that came about 8 years after its discovery, after \cite{Jones-Annals} appeared in print.   Following  ideas pioneered by V.I. Arnold, and carried out by his then-student V. Vassiliev\cite{Vassiliev}, let's now consider the collection of {\it all} knots in $\mathbb S^3$.  Vassiliev found a way to give mathematical precision to that concept.    His first change was a very natural one in the history of mathematics: instead of thinking of a knot as the image  of $S^1$ under an embedding $e:\mathbb S^1\to\mathbb R^3$, he focused on the embedding itself.  That allowed him to broaden his viewpoint and look at the space $\mathcal M$ of all smooth {\it maps} from $\mathbb S^1$ to $\mathbb R^3$.   Define the {\it discriminant} $\Sigma$ of $\mathcal M$ to be the subset of maps that are {\it not} embeddings.  Vassiliev studied the components of $ \mathcal M \setminus \Sigma$, a set that has a component for each knot type.  See the excellent introductory material in \cite{Mortier} for appropriate  restrictions on admissible maps, and topologies on the spaces themselves,  allowing him to study $\mathcal M \setminus \Sigma$  using known tools, and for an excellent new mathematical contribution to a major area of mathematics.

Vassiliev's invariants, rooted in the {\it Cohomology} of $\mathcal M \setminus \Sigma$, have an {\it order} that can be described in the following intuitive way:  Any knot can be changed to the unknot by a series of transverse crossing changes,  and so  any knot  can be changed to any other by crossing changes.  The space $\mathcal M$ is therefore connected, and the minimum number of crossing changes serves as a measure of complexity, the {\it order} of  an element of $\mathcal M \setminus \Sigma$.  On the way, one looks at {\it singular knots}, that is knots that have a representative with one or more transverse double points,  interesting mathematical objects in their own right. The following theorem was proved in \cite{Birman-Lin}:
\begin{theorem} \label{Th:J detects V} Let $\bf K$ be a knot in $\mathbb R^3$, and let $V_{\bf K(t)}$ be its 1-variable Jones polynomial.  Let  $U_{\bf K}$ be the power series obtained from $V_{\bf K(t)}$ by replacing the variable $t$ in the Jones polynomial by $e^x$, and then rewriting $U_{\bf K}$ as a power series in $x$: 
$$U(\bf K) = \sum_{i=0}^{\infty}u_i(\bf K) x^i$$
Then $u_0(\bf K) = 1$ and each $u_i(\bf K), i >1$, is a Vassiliev invariant of order $i$.
\end{theorem}
\begin{openproblem} Is there a structure in the algebra part of the dictionary that is related to Jones' work and 
corresponds to  Vassiliev's space of all knots?
\end{openproblem}

\section{Back to algebra} \label{s:quantum topology}  Returning to the Jones polynomial $V_{\bf K}(t)$ of a knot type $\bf K$ in $\mathbb R^3$, I was reminded, as I was writing this paper, of the explosion of ideas that followed. To bring a semblance of order out of that particular chaos, and focus on the part of the story that I wish to tell next,  recalling the following sequence of discoveries that occurred after May 31,1984, when the only act in town  was $V_{\bf K}(t)$:   
\begin{enumerate}
\item A 2-variable polynomial invariant $P_{\bf K}(t,x)$ of a knot ${\bf K}$ in 3-space appeared on the scene.  See $\S$6, and $\S$11-13 of \cite{Jones-Annals}.  The associated algebra, in the 2-variable case, is the Hecke algebra 
$H_n(l,m)$ of type $A_{n-1}$. Its irreducible representations are in 1-1 correspondence with the Young diagrams that appear in the literature about the symmetric group.  The Temperley-Lieb algebra,  that is the algebra associated to all Young diagrams with exactly 2 rows, gave rise to the   1-variable polynomial $V_{\bf K}(t)$ that was Jones' initial discovery. 
\item Louis Kauffman gave a a very simple new proof of the existence of $V_{\bf K}(t)$ in \cite{Kauffman-Topology}, using states-models techniques. His methods are 
combinatorial in nature, and related to von Neumann algebras in the sense that Kauffman gave a diagrammatic (planar algebra) interpretation of the role of the Temperley-Lieb algebra in the construction of the Jones polynomial.
\item Having in hand $P_{\bf K}(t,x)$ and states-models techniques, Kauffman went on to discover yet another 2-variable knot polynomial, now known as  the {\it Kauffman polynomial} $K_{\bf K}(l,m)$. 
For this I refer the reader to Kauffman's article  about its discovery \cite{Kauffman-JKTR}. That brings me back to Vaughan Jones' letter-writing habits, because the reader will find at the website {\it celebratio.org} a long letter from Jones to Kauffman about the Kauffman polynomial, and a discussion of proposed ramifications in Physics.  
\end{enumerate}

Knowing about the Homfly and Kauffman polynomials,  and with knowledge of the role of the Hecke algebra in the situation of  the Homfly polynomial, it was natural that Hans Wenzl (a graduate student working with Jones, when he and I first met) and I, chatting, speculated that perhaps the Kauffman polynomial $K_{\bf K}(l,m)$ could itself be derived from an appropriate algebra.  In other words, was there a new word in the algebra side of the dictionary that was associated to the newly discovered Kauffman polynomial on the topology side?    The answer is `yes', and the new algebra is known as the {\it BMW Algebra}, where the `M' stands for the fact that some of the results in \cite{Birman-Wenzl} (but not all) were discovered simultaneously and independently by J. Murakami \cite{Murakami}.   An unexpected surprise was that the new algebra turned out to be related to a known one, that is the centralizer algebra constructed by Brauer in the 1930's. In this way the interconnections between Algebra and Topology, the theme of this note, were expanded and deepened.  

We end by mentioning one more open problem that relates to Jones' discoveries:
\begin{openproblem}  A major result in  {\rm \cite{Jones-Index}}, which started the topic of this note, is that the index of a  subfactor of a type {\rm II}$_1$ factor is either $\geq 4$ or 
${4  cos^2 (\pi/n) = \{1,2, 3/2+\sqrt{5}/2, 3, \dots} \}$ for some $n = 3,4,5,\dots$.  Moreover, all of these values occur. 
What is the meaning of these restrictions, if any,  on the topology side of the dictionary?  
\end{openproblem}

Joan S. Birman\\
jb@math.columbia.edu\\
Prof. Emeritus, Math Dept., Barnard College\\
11 Riverside Drive, Apt 3BW\\
New York, NY 10023\\

\end{document}